\documentclass[12pt,leqno]{article}

\usepackage{amsfonts}
\usepackage{amsmath}
\usepackage{amssymb}
\usepackage[dvips]{graphicx}

\newcommand{\Xcomment}[1]{}

\newtheorem{theorem}{Theorem}[section]
\newtheorem{lemma}[theorem]{Lemma}
\newtheorem{corollary}[theorem]{Corollary}

\newenvironment{numitem}{\refstepcounter{equation}\begin{enumerate}%
\item[(\arabic{equation})]$\quad$}{\end{enumerate}}

\newenvironment{proof}{\noindent {\bf Proof.\/}}{$\qed$\bigskip}

\newcommand{\refeq}[1]{(\ref{eq:#1})}  

\def\qed{ \ \vrule width.2cm height.2cm depth0cm}

\def\tilde{\widetilde}
\def\hat{\widehat}
\def\bar{\overline}
\def\eps{\epsilon}

\def\Rset{{\mathbb R}}
\def\Zset{{\mathbb Z}}

\def\Cscr{\mathcal{C}}

\def\Fscr{\mathcal{F}}

\def\Pscr{\mathcal{P}}
\def\Qscr{\mathcal{Q}}
\def\Tscr{\mathcal{T}}

\def\diver{{\rm div}}

\def\deltain{\delta^{\rm in}}
\def\deltaout{\delta^{\rm out}}
\def\degin{{\rm deg}^{\rm in}}
\def\degout{{\rm deg}^{\rm out}}
\def\valu{{\rm val}}
\def\supp{{\rm supp}}

\begin{document}

\begin{titlepage}

\title{\Large Free multiflows in bidirected and
skew-symmetric graphs}

\author{
Maxim~A.~Babenko\thanks{Dept. of Mechanics and Mathematics,
Moscow State University, Vorob'yovy Gory, 119899 Moscow,
Russia, {\sl email}: mab@shade.msu.ru.
Supported by RFBR grants 03-01-00475 and NSh 358.2003.1.}
  \and
Alexander~V.~Karzanov\thanks{
Institute for System Analysis, 9, Prospect 60 Let Oktyabrya,
117312 Moscow, Russia, {\sl email}: sasha@cs.isa.ru. Supported by
NWO--RFBR grant 047.011.2004.017.}
}
\maketitle

\begin{abstract}
A graph (digraph) $G=(V,E)$ with a set $T\subseteq V$
of terminals is called {\em inner Eulerian} if each nonterminal
node $v$ has even degree (resp. the numbers of edges entering and
leaving $v$ are equal). Cherkassky~\cite{cher-77} and
Lov\'asz~\cite{lov-76} showed that the maximum number of pairwise
edge-disjoint $T$-paths in an inner Eulerian graph $G$ is equal to
$\frac12\sum_{s\in T} \lambda(s)$, where
$\lambda(s)$ is the minimum
number of edges whose removal disconnects $s$ and $T-\{s\}$. A
similar relation for inner Eulerian digraphs was established by
Lomonosov~\cite{lom-85}. Considering undirected and directed
networks with ``inner Eulerian'' edge
capacities, Ibaraki, Karzanov, and Nagamochi~\cite{IKN-98} showed
that the problem  of finding a maximum integer multiflow (where
partial flows connect arbitrary pairs of distinct terminals) is
reduced to $O(\log T)$ maximum flow computations and to a number of
flow decompositions.

In this paper we extend the above max-min relation to inner Eulerian
{\em bidirected} and {\em skew-symmetric} graphs and develop
an algorithm of complexity $O(VE\log T\log(2+V^2/E))$ for the
corresponding capacitated cases. In particular, this improves the
bound in~\cite{IKN-98} for digraphs. Our algorithm uses a fast
procedure for decomposing a flow with $O(1)$ sources and sinks in a
digraph into the sum of one-source-one-sink flows.
\end{abstract}

\medskip
\noindent {\em KeyWords}:
bidirected graph, skew-symmetric graph, edge-disjoint paths,
multiflow.

\medskip
\noindent {\em AMS Subject Classification}:
90C27, 90B10

\def\thepage {} 

  \end{titlepage}


\section{{\Large\rm Introduction}} \label{sec:intr}

A graph (digraph) $G=(V,E)$ with a distinguished subset $T$ of
nodes is said to be {\em inner Eulerian} if each node $v\in V-T$
has even degree (resp. the indegree and outdegree of $v$ are
equal). The nodes in $T$ and in $V-T$ are called {\em terminals}
and {\em inner} nodes, respectively. A simple path in $G$ is called
a $T$-{\em path} if its ends are distinct terminals and the other
nodes are inner. There is a nice max-min relation established by
Cherkassky~\cite{cher-77} and Lov\'asz~\cite{lov-76} for graphs,
and by Lomonosov~\cite{lom-85} for digraphs:
  \begin{numitem}
if $(G,T)$ is inner Eulerian, then the maximum number $\nu_{G,T}$ of
pairwise edge-disjoint $T$-paths is equal to
$\frac12\sum_{s\in T}\lambda_{G,T}(s)$.
  \label{eq:CLL}
  \end{numitem}

Here and later on, for a subset $X$ of nodes,
$\delta(X)=\delta_G(X)$ denotes the set of edges with one end in
$X$ and the other in $V-X$, called the {\em cut} induced by $X$. For
$s\in T$, we refer to a subset $X\subset V$ with $X\cap T=\{s\}$ as
an $s$-{\em set}. Then $\lambda_{G,T}(s)$ is defined to be the
minimum cardinality $|\delta(X)|$ among the $s$-sets $X$.

The above max-min relation has an obvious extension to the
capacitated case. Given a nonnegative integer function $c:E\to\Zset_+$
of edge {\em capacities},
let us say that the triple ({\em network}) $(G,T,c)$
is {\em inner Eulerian}
if for each inner node $v$, the total capacity of
edges incident with $v$ is even when $G$ is a graph, and the total
capacity of edges entering $v$ is equal to that of edges leaving $v$
when $G$ is a digraph. Then~\refeq{CLL} yields the following
relation for an inner Eulerian $(G,T,c)$:
  \begin{equation} \label{eq:mm}
\max\{\valu(\Fscr)\}=\frac12\sum\nolimits_{s\in T}\lambda_{c,T}(s),
  \end{equation}
where $\lambda_{c,T}(s)$ denotes the minimum cut capacity
$c(\delta(X))$ among the $s$-sets $X$, the maximum is taken
over all collections $\Fscr$ of $T$-paths $P_1,\ldots,P_k$ along
with nonnegative {\em integer} weights $\alpha_1,\ldots,\alpha_k
\in\Zset_+$ that satisfy the packing condition
  \begin{equation}  \label{eq:pack}
\sum(\alpha_i: e\in P_i)\le c(e) \qquad\mbox{for all $e\in E$},
  \end{equation}
and $\valu(\Fscr)$ denotes the {\em total value} of $\Fscr$, defined to
be $\alpha_1+\ldots+\alpha_k$%
\footnote{Originally relation~\refeq{mm} was stated for fractional
multiflows $\Fscr$ in~\cite{kup-71}, with a flaw in the proof.
  }.
(Hereinafter for a function $g:S\to\Rset$ and a subset
$S'\subseteq S$, $g(S')$ stands for $\sum_{e\in S'}g(e)$.)

A collection $\Fscr$ consisting of $T$-paths $P_i$ with real weights
$\alpha_i\in\Rset_+$ that obeys~\refeq{pack} is called a {\em free
multiflow}  (the adjective ``free'' is used to emphasize that any pair
of distinct terminals is allowed to be connected by a path,
i.e., the commodity graph in the multiflow maximization problem
is complete).
A multiflow achieving the equality in~\refeq{mm} is called
{\em maximum}. Thus, whenever $(G,T,c)$ is inner Eulerian, there
exists an {\em integer} maximum free multiflow (i.e., having the weights
of all paths integral).

Cherkassky~\cite{cher-77} showed that such a
multiflow in an inner Eulerian undirected network can be found
in strongly polynomial time. Subsequently much faster algorithms
both for graphs and digraphs have been developed. They apply
a ``divide-and-conquer'' approach in which a current
network $(G,T,c)$ with $|T|\ge 4$ is recursively replaced by two
networks $(G',T',c')$ and $(G'',T'',c'')$ such
that $|T'|,|T''|\le\lceil |T|/2\rceil+1$.
Originally, such an approach was applied in~\cite{kar-79}
to find, in $O(\phi(V,E)\log T)$ time, a {\em half-integer} maximum
multiflow in a graph $G$ with integer edge capacities (but
not guaranteeing integrality in the inner Eulerian case).
Hereinafter, in notation involving functions of numerical arguments
or time bounds, we indicate sets for their cardinalities, and
$\phi(n,m)$ stands for the complexity of an algorithm for finding a
maximum flow in a network with $n$ nodes and $m$ edges.

This algorithm was improved and extended in~\cite{IKN-98} so as to
find an integer maximum free multiflow in an inner Eulerian
undirected network in the same time $O(\phi(V,E)\log T)$, and in an
inner Eulerian directed network in $O(\phi(V,E)\log T+V^2E)$ time.

\medskip
\noindent{\bf Remark 1.} The inner Eulerianness condition is
important. Withdrawing it makes the undirected problem more
difficult, though still polynomially solvable in the noncapacitated
case (a max-min relation is due to Mader~\cite{mad-78} and an
original polynomial algorithm is due to Lov\'asz~\cite{lov-80}),
and makes the directed noncapacitated problem NP-hard already
for two terminals~\cite{FHW-80}.

\medskip
The purpose of this paper is to extend the above theoretical and
algorithmic results to bidirected graphs. (This sort of nonstandard
graphs was introduced by Edmonds and Johnson~\cite{EJ-70} in
connection with one important class of integer linear programs
generalizing problems on flows and matchings; for
a survey, see also~\cite{law-76,sch-03}.)

\begin{figure}[tb]
\begin{center}
\includegraphics{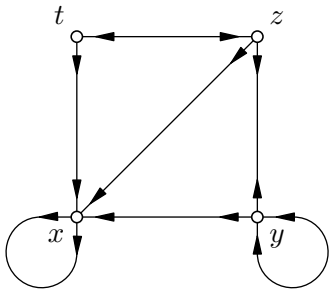}
\end{center}
\caption{A bidirected graph instance}
\label{fig:bidir}
\end{figure}

Recall that in a {\em
bidirected} graph $G=(V,E)$ three types of edges are allowed:
(i) a usual directed edge, or an {\em arc}, that leaves one node
and enters another one; (ii) an edge {\em from both} of its ends;
or (iii) an edge {\em to both} of its ends.
When $u=v$, the edge becomes a loop; in what follows we
admit only loops of types (ii) and (iii) (as loops of type (i) do
not affect our problem and can be excluded from consideration).
A nonloop edge entering a node $v$ contributes one unit to the
indegree $\degin(v)$ of $v$, while a loop of type (iii) at $v$
contributes two units to $\degin(v)$; the outdegree $\degout(v)$ of
$v$ is specified in a similar way.
Edges $e,e'$ connecting nodes $u,v$ are called {\em parallel} if
$e$ enters $u$ if and only if $e'$ does so, and similarly for $v$.
If $G$ has no parallel edges, then $|E|\le 2|V|^2$.
An instance of bidirected graphs is
drawn in Fig.~\ref{fig:bidir}.

The notion of inner Eulerianness for a bidirected graph $G$
with a set $T$ of terminals is analogous to that for usual
digraphs: $\degin(v)=\degout(v)$ for
all inner nodes $v$. Inner Eulerian triples $(G,T,c)$, where
$c:E\to\Zset_+$, are those that turn into inner Eulerian pairs
$(G',T)$ when each edge $e$ is replaced by $c(e)$ parallel edges.

In order to be able to extend the above results to bidirected
graphs, we need to admit $T$-paths with restricted self-intersections.
(In undirected or directed graphs, when one node is reachable
from another one by a path, then it is reachable by a simple path
as well, but this need not hold in a bidirected graph.)
A {\em walk} in a bidirected graph $G$ is an alternating sequence
$P=(s=v_0,e_1,v_1,\ldots,e_k,v_k=t)$ of nodes and edges such that
each edge $e_i$ connects nodes $v_{i-1}$ and $v_i$, and for
$i=1,\ldots,k-1$, the edges $e_i,e_{i+1}$ form a {\em transit pair}
at $v_i$, which means that one of $e_i,e_{i+1}$ enters and the
other leaves $v_i$. Note that $e_1$ may enter $s$ and $e_k$ may
leave $t$; nevertheless, we refer to $P$ as a walk from $s$ to $t$,
or an $s$--$t$ {\em walk}. $P$ is a {\em cycle} if $v_0=v_k$ and
the pair $e_1,e_k$ is transit at $v_0$; a cycle is usually
considered up to cyclically shifting. Observe that an $s$--$s$ walk
is not necessarily a cycle. By a {\em path} in a bidirected graph $G$
we mean an {\em edge-simple} walk $P$, i.e., a walk with all
edges different. Similar to usual graphs/digraphs,
a $T$-path (a $T$-walk) is meant to be a path (resp. walk) whose
ends are distinct terminals and the other nodes are inner.

Define the number $\lambda_{G,T}(s)$ for $s\in T$ as before.
We show the following

  \begin{theorem} \label{tm:main}
Property \refeq{CLL} remains valid for a bidirected graph $G$ and
the set of $T$-paths in it.
  \end{theorem}

\noindent{\bf Remark 2.} In this theorem it suffices to consider
only minimal $T$-paths, where a path (edge-simple cycle)
$P=(v_0,e_1,v_1,\ldots,e_k,v_k)$ is
called {\em minimal} if no part of $P$ from $v_i$ to $v_j$ with
$i<j$ (resp. $0<j-i<k$) forms a cycle. (A minimal path/cycle need
not be simple but it passes any node of $G$ at most twice.)
Moreover, one can consider only those $T$-paths whose induced
bidirected graphs contain no cycle. (In the underlying undirected
graph of such a path, each edge belongs to at most one circuit.)

\medskip
A usual digraph is a special case of bidirected graphs and
Theorem~\ref{tm:main} generalizes the above-mentioned result
in~\cite{lom-85}. Also there is a natural correspondence between
the $T$-paths in an undirected graph $G$ and the
minimal $T$-paths in the bidirected graph $G'$ formed from $G$ as
follows: direct each edge of $G$ from both of its ends,
and for each inner node $v$, assign $\lceil\deg(v)/2\rceil$
loops entering (twice) $v$. Then a $T$-path $P$ in $G$
is turned into a $T$-path in $G'$ by adding one loop to each
intermediate node of $P$. Moreover, $(G',T)$ is inner Eulerian
if $(G,T)$ is such. Due to this correspondence,
Theorem~\ref{tm:main} generalizes the above-mentioned
Cherkassky--Lov\'asz' result for undirected graphs as well.

Like the pure graph and digraph cases, one can reformulate
Theorem~\ref{tm:main} in capacitated terms: relation~\refeq{mm}
concerning integer free multiflows $\Fscr$ remains valid when $G$ is
bidirected and $(G,T,c)$ is inner Eulerian. In this case one should
consider $T$-walks, rather than $T$-paths, and refine the packing
condition~\refeq{pack} as
  \begin{equation} \label{eq:packbd}
\sum_{i=1}^k \alpha_in_i(e)\le c(e) \qquad\mbox{for all $e\in E$},
  \end{equation}
where $n_i(e)$ is the number of occurrences of an edge $e$ in a walk
$P_i$. Thus, the above problem for undirected and directed networks
is generalized as:
  \begin{itemize}
\item[(P)]
{\em Given an inner Eulerian network $(G,T,c)$, where $G$ is
bidirected, find a collection (free multiflow) $\Fscr$ of
$T$-walks $P_1,\ldots,P_k$ with weights
$\alpha_1,\ldots,\alpha_k\in\Zset_+$ that satisfies~\refeq{packbd}
and maximizes the value $\valu(\Fscr):=\alpha_1+\ldots+\alpha_k$.}
  \end{itemize}

\noindent{\bf Remark 3.} Let $X$ be an arbitrary subset of nodes of $G$.
One can modify $G$ as follows: for each node $v\in X$ and each edge
$e$ incident with $v$, reverse the direction of $e$ at $v$. Also
for an arbitrary arc $e$ incident with a terminal $s$, one can
reverse the direction of $e$ at $s$. Both transformations preserve
the inner Euleriannes of $(G,T,c)$ and the set of $T$-walks.
Therefore, problem~(P) remains equivalent under such
transformations.

\medskip
Another appealing class of nonstandard graphs was introduced by
Tutte~\cite{tut-67} who originated a mini-theory, parallel
to~\cite{EJ-70} in a sense, aiming to unify and generalize flow and
matching problems. These are so-called {\em skew-symmetric graphs}
(or anti-symmetrical digraphs, in Tutte's terminology), digraphs with
involutions on the nodes and on the arcs which reverse the
orientation of each arc (a precise definition is given in
Section~\ref{sec:skew}). His and other researchers' study of
structural and optimization problems on skew-symmetric graphs has
resulted in a number of interesting theorems, methods and
applications.

There is a close relationship
between skew-symmetric and bidirected graphs, and typically results
on the former can be reformulated for the latter, and vice versa.
So is for the problem of our study, too. We take advantage from
both representations. The language of bidirected graphs is more
preferable for us to work in the non-capacitated case; we prove
Theorem~\ref{tm:main} directly and obtain its analog for
skew-symmetric graphs as a corollary. On the other hand, we prefer
to deal with skew-symmetric graphs in algorithmic design for the
capacitated case. (A serious reason is that a flow in a bidirected
network is defined as a packing of $T$-paths and we do not see
reasonable alternative settings for it, while a flow in a
skew-symmetric network can be given in a more compact form, via a
function on the arc set.) Some facts about skew-symmetric
flows and technical tools elaborated for them help us to
devise a fast algorithm for the skew-symmetric analog of
problem~(P) concerning {\em integer skew-symmetric free multiflows}
in an inner Eulerian skew-symmetric network. This yields a fast
algorithm for~(P) as well.

The paper is organized as follows. In Section~\ref{sec:proof} we
prove Theorem~\ref{tm:main} (which is relatively simple) relying
on the fact that an inner Eulerian bidirected graph can be
decomposed into cycles and paths with both ends in $T$.
Section~\ref{sec:skew} explains the correspondence between bidirected
and skew-symmetric graphs, reviews some known results about
the latter (in particular, Tutte's result on symmetric
decompositions of skew-symmetric flows) and gives a skew-symmetric
analog of Theorem~\ref{tm:main}.
Section~\ref{sec:alg} develops an algorithm
for finding a maximum integer skew-symmetric free multiflow in an
inner Eulerian skew-symmetric network. It relies on a general approach
in~\cite{kar-79} and some ingredients from~\cite{IKN-98} and attracts
additional combinatorial ideas and techniques. As a consequence,
problem~(P) is solved in time $O(VE\log T\log(2+V^2/E))$ (if
the $O(nm\log(2+n^2/m))$-algorithm of Goldberg and Tarjan~\cite{GT-88}
is applied for finding a maximum flow in a directed network with
$n$ nodes and $m$ arcs).
This improves the bound in~\cite{IKN-98} for digraphs.
To achieve this bound, we use a faster procedure for the particular
flow decomposition problem: given an integer flow $f$ with $O(1)$
sources and sinks in a digraph with $n$ nodes and $m$ arcs,
decompose $f$ into the sum of integer flows, each connecting one
source to one sink. The procedure developed in
Section~\ref{sec:decom} solves this problem in $O(m\log(2+n^2/m))$
time. In the concluding Section~\ref{sec:sdecom}, this procedure is
extended to symmetric flows in skew-symmetric graphs (it is not
used in the algorithm for~(P) but may be of interest for other
applications).

\section{{\Large\rm Proof of Theorem~\ref{tm:main}}}
\label{sec:proof}

Let $G=(V,E)$ be an inner Eulerian bidirected graph with
a terminal set $T$.
One may assume that $G$ has no loops incident with terminals.
Since $(G,T)$ is inner Eulerian, for each inner node $v$, one can
choose a set $\pi(v)$ of transit pairs at $v$ so that each non-loop
edge incident with $v$ occurs in exactly one pair and each loop at
$v$ (if any) occurs in two pairs.
The collection $\{\pi(v):v\in V-T\}$  determines a
decomposition of (the edge-set of) $G$ into a collection $\Cscr$ of
edge-simple cycles and a collection $\Pscr$ of paths with both ends
at $T$. More precisely, each edge $e\in E$ belongs to exactly one
member $P$ of $\Cscr\cup\Pscr$ and satisfies the following condition:
for each end $v$ of $e$, if $v\in V-T$ and $\{e,e'\}\in\pi(v)$, then
either $e,v,e'$ or $e',v,e$ are three consecutive elements in $P$,
while if $v\in T$, then $P$ begins with $v,e$ or ends with $e,v$.
Note that all nodes of any cycle in $\Cscr$
and all intermediate nodes of any path in $\Pscr$ are inner. So each
path in $\Pscr$ is a $T$-path unless it connects equal terminals.
When needed, we may reverse some paths in $\Pscr$.

For $s\in T$, let $\Pscr_s$ ($\Qscr_s$) denote the set of paths in
$\Pscr$ with exactly one end (resp. with both ends) at $s$.
Since $|\Pscr_s|+2|\Qscr_s|=\deg(s)
\ge\lambda_{G,T}(s)$ (where $\deg(v)$ is the full
degree $\degout(v)+\degin(v)$ of $v$), the theorem becomes
trivial when all sets $\Qscr_s$ are empty. In a general case, we
try to transform the decomposition so as to increase the ``useful
value'' $\eta(\Pscr,\Cscr):=\sum_{s\in T}|\Pscr_s|$, by applying a
certain augmenting approach.

Consider $s\in T$ and assume, w.l.o.g., that all paths in $\Pscr_s$
begin at $s$. Let $L=(x_0,x_1,\ldots,x_q)$ be a sequence of distinct
nodes such that
  \begin{numitem}
either $L=\{s\}$, or $x_0$ belongs to a path in $\Qscr_s$,
and for $i=1,\ldots,q$, the nodes $x_i,x_{i-1}$ occur in a cycle in
$\Cscr$ or occur in {\em this order} in a path in $\Pscr_s$.
  \label{eq:aug}
  \end{numitem}
We say that $L$ is augmenting if $x_q$ belongs to a path in $\Pscr$
having both ends in $T-\{s\}$. Consider two cases.

\medskip
{\em Case 1}. There is no augmenting sequence for $s$. Let $X_s$ be
the set of all nodes occurring in sequences $L$ as in~\refeq{aug}.
Clearly $X_s\cap T=\{s\}$.
Consider an edge $e$ of the cut $\delta(X_s)$; let $u,v$ be the
ends of $e$ in $X_s$ and $V-X_s$, respectively.
Observe that $e$ belongs to neither a cycle in
$\Cscr$ nor a path in $\Pscr-\Pscr_s$. Also (by~\refeq{aug}) if $e$
belongs to a path $P\in\Pscr_s$, then all nodes of $P$ from $s$ to
(the last occurrence of) $u$ are contained in $X_s$, i.e., $P$
traverses the cut $\delta(X_s)$ exactly once. This implies
$|\Pscr_s|=|\delta(X_s)|$.

Hence, if none of terminals admits an augmenting sequence as above,
then the number of $T$-paths in $\Pscr$ is at least
$\frac12\sum_{s\in T}\lambda_{G,T}(s)$, as required.

\medskip
{\em Case 2.} An augmenting sequence $L=(x_0,x_1,\ldots,x_q)$ for
$s$ exists. Let $L$ be chosen so that no proper subsequence in it
is augmenting. Then:
  \begin{numitem}
any cycle in $\Cscr$ meets at most two nodes in $L$ and these
nodes are consecutive in $L$;
  \label{eq:intC}
  \end{numitem}
  \begin{numitem}
if a path $P\in\Pscr_s$ contains a node $x_i$, then the part of $P$
from $s$ to (the last occurrence of) $x_i$ can contain at most one
node $x_j$ with $j>i$; moreover, if such an $x_j$ exists then
$j=i+1$.
  \label{eq:intP}
  \end{numitem}
We transform $(\Pscr,\Cscr)$ along $L$, step by step, as follows.
Choose a $Q\in\Qscr_s$ containing $x_0$.
At the first step, if (a) $x_0,x_1$ belong to a cycle $C\in\Cscr$,
then we combine $Q$ and $C$ into one $s$--$s$ path.
And if (b) $x_0,x_1$ belong to a path $P\in\Pscr_s$ from $s$ to $t$,
say, and if $x_1$ occurs in $P$ earlier than $x_0$, then we replace
$Q$ by the concatenation of the part $P'$ of $P$ from $s$ to (the last
occurrence of) $x_0$ and the part $Q'$ of $Q$ from $x_0$ to $s$, and
replace $P$ by the concatenation of the rest of $Q$ (from $s$ to
$x_0$) and the rest of $P$ (from $x_0$ to $t$). (We assume, w.l.o.g.,
that the last edge of $P'$ and the first edge of $Q'$ form a
transit pair at $x_0$; otherwise reverse $Q$.) As a result, we
obtain an $s$--$s$ path, denoted by $Q$ as before, that contains
$x_1$. In case (a), the cycle $C$ vanishes, and in case (b), the
new path $P$ goes from $s$ to $t$ as before, and its part
from $x_0$ to $t$ preserves.
This together with~\refeq{intC} and~\refeq{intP} implies validity
of~\refeq{aug} for the remaining sequence $(x_1,\ldots,x_q)$;
moreover, \refeq{intC} and~\refeq{intP} are maintained as well. At
the second step, we consider the pair $x_1,x_2$ and act in a
similar way, and so on.

Eventually, after $q$ steps, the current $s$--$s$ path $Q$ contains
the node $x_q$. Since $L$ is augmenting, $x_q$ also belongs to some
$t$--$p$ path $R\in\Pscr$ with $t,p\in T-\{s\}$ (possibly $p=t$).
Now splitting $Q,R$ at $x_q$ and concatenating the arising four
pieces in another way, we obtain two $T$-paths, one connecting $s$
and $p$ and the other connecting $s$ and $t$. This gives a new
decomposition $(\Pscr,\Cscr)$ of $G$ having a larger value of $\eta$,
and the theorem follows. \qed

\medskip
In fact, the above proof is constructive and prompts a polynomial
algorithm for finding a maximum number of pairwise edge-disjoint
$T$-paths in an inner Eulerian bidirected graph. A more efficient
and more general algorithm (dealing with the capacitated case) is
described in Section~\ref{sec:alg}.

\section{{\Large\rm Skew-Symmetric Graphs}} \label{sec:skew}

This section contains terminology and some basic facts concerning
skew-symmetric graphs and explains the correspondence between these and
bidirected graphs. For a more detailed survey on skew-symmetric graphs,
see, e.g., ~\cite{tut-67,GK-96,GK-04}.

A {\em skew-symmetric graph} is a digraph $G=(V,E)$ endowed with
two bijections $\sigma_V,\sigma_E$ such that: $\sigma_V$ is
an involution on the nodes (i.e., $\sigma_V(v)\ne v$ and
$\sigma_V(\sigma_V(v))=v$ for each $v\in V$), $\sigma_E$ is an
involution on the arcs, and for each arc $e$ from $u$ to $v$,
$\sigma_E(e)$ is an arc from $\sigma_V(v)$ to $\sigma_V(u)$. For
brevity, we combine the mappings $\sigma_V,\sigma_E$ into one mapping
$\sigma$ on $V\cup E$ and call $\sigma$ the {\em symmetry} (rather
than skew-symmetry) of $G$.
For a node (arc) $x$, its symmetric node (arc) $\sigma(x)$ is also
called the {\em mate} of $x$, and we will often use notation with
primes for mates, denoting $\sigma(x)$ by $x'$. Obviously,
$\degin(v)=\degout(v')$ for each $v\in V$.

We admit parallel arcs, but not loops, in $G$.
Observe that if $G$ contains an arc $e$ from a node $v$
to its mate $v'$, then $e'$ is also an arc from $v$ to $v'$ (so the
number of arcs of $G$ from $v$ to $v'$ is even and these parallel
arcs are partitioned into pairs of mates).

By a path (circuit) in $G$ we mean a
simple directed path (cycle), unless explicitly stated otherwise.
The symmetry $\sigma$ is extended in a natural way to paths,
subgraphs, and other objects in $G$. In particular, two paths or
circuits are symmetric to each other if the elements of one of them
are symmetric to those of the other and go in the reverse order:
for a path (circuit) $P=(v_0,e_1,v_1,\ldots,e_k,v_k)$, the symmetric
path (circuit) $\sigma(P)$ is $(v'_k,e'_k,v'_{k-1},\ldots,e'_1,v'_0)$.
One easily shows that $G$ cannot contain self-symmetric circuits
(cf.~\cite{GK-04}).

Following terminology in~\cite{GK-96}, a path or circuit in $G$ is
called {\em regular} if it contains no pair of symmetric arcs
(while symmetric nodes in it are allowed). For a function $h$ on
$E$, its symmetric function $h'$ is defined by $h'(e'):=h(e)$,
$e\in E$, and $h$ is called (self-){\em symmetric} if $h=h'$.

For a function $f:E\to \Rset$ and a node $v\in V$, define
  $$
\diver_f(v):=\sum(f(e):e\in\deltaout(v))-\sum(f(e):e\in\deltain(v)),
  $$
(the {\em divergency} of $f$ at $v$),
where $\deltaout(v)$ ($\deltain(v)$) denotes the set of arcs of $G$
leaving (resp. entering) $v$. Let $f$ be nonnegative, integer-valued
and symmetric, and let $S$ be a subset of nodes not intersecting
$S'=\sigma(S)$. When $\diver_f(v)$ is nonnegative at each $v\in S$
and zero at each $v\in V-(S\cup S')$, $f$ is said to be an
{\em IS-flow} (integer symmetric flow) from $S$ to $S'$. The
{\em value} $\valu(f)$ of $f$ is $\sum_{s\in S}\diver_f(s)$.
By a multiterminal version of a theorem due to
Tutte~\cite{tut-67}, an IS-flow $f$ from $S$ to $S'$ has an
{\em integer symmetric decomposition}. This means that
  \begin{numitem}
$f$ is representable as $f=\alpha_1\chi^{P_1}+\alpha_1\chi^{P'_1}+
\ldots+\alpha_k\chi^{P_k}+\alpha_k\chi^{P'_k}$, where for
$i=1,\ldots,k$, $P_i$ is a path from $S$ to $S'$ or a circuit,
$P'_i$ is the path (also going from $S$ to $S'$) or circuit symmetric
to $P_i$, and $\alpha_i\in\Zset_+$.
  \label{eq:sdecom}
  \end{numitem}
Here $\chi^P$ denotes the incidence vector of the arc-set
of a path/circuit $P$, i.e., for $e\in E$, $\chi^P(e)=1$ if $e$
belongs to $P$, and 0 otherwise.
Note that paths/circuits in~\refeq{sdecom} need not be regular.
Considering $\alpha_i$ as the weight of $P_i$ and of $P'_i$, observe
that the total weight of paths from $S$ to $S'$ is equal to
$\valu(f)$. Similar to flow decomposition in usual digraphs, an
integer symmetric decomposition of an IS-flow $f$ can be found in
$O(VE)$ time.

Let $(G,T)$ be inner Eulerian, where the terminal set $T$ is
(self-)symmetric. Take a partition $\{S,S'=\sigma(S)\}$ of $T$
such that $\degout(s)\ge\degin(s)$ for all $s\in S$. Since the
all-unit function $f$ on $E$ represents an IS-slow from
$S$ to $S'$, \refeq{sdecom} implies that
  \begin{numitem}
there exists a symmetric collection $\Pscr$ of circuits and paths
from $S$ to $S'$ in $G$ such that the members of $\Pscr$ are pairwise
arc-disjoint and cover $E$, and each terminal $s\in S$ is
the beginning of exactly $\degout(s)-\degin(s)$ paths in $\Pscr$.
  \label{eq:decomp}
  \end{numitem}
Moreover, the members of $\Pscr$ are regular (for if some
$P\in\Pscr$ contains mates $e,e'\in E$, then $e,e'$ are in
$\sigma(P)$ as well, which is impossible).

\medskip
Next we explain the correspondence between skew-symmetric and
bidirected graphs (cf.~\cite[Sec.~2]{GK-04}). For sets $X,A,B$, we
may use notation $X=A\sqcup B$ when $X=A\cup B$ and $A\cap B=
\emptyset$. Given a skew-symmetric
graph $G=(V,E)$, choose an arbitrary partition $\pi=\{V_1,V_2\}$ of
$V$ such that $V_2$ is symmetric to $V_1$. Then $G,\pi$ determine
bidirected graph $H$ with node set $V_1$ whose edges
correspond to the pairs of symmetric arcs in $G$. More precisely,
arc mates $a,a'$ of $G$ generate one edge $e$ of $H$ connecting nodes
$u,v\in V_1$ such that: (i) $e$ goes from $u$ to $v$ if one of $a,a'$
goes from $u$ to $v$ (and the other goes from $v'$ to $u'$ in
$V_2$); (ii) $e$ leaves both $u,v$ if one of $a,a'$ goes from $u$
to $v'$ (and the other from $v$ to $u'$); (iii) $e$ enters both
$u,v$ if one of $a,a'$ goes from $u'$ to $v$ (and the other from
$v'$ to $u$). In particular, $e$ is a loop if $a,a'$ connect a pair
of symmetric nodes.

\begin{figure}[tb]
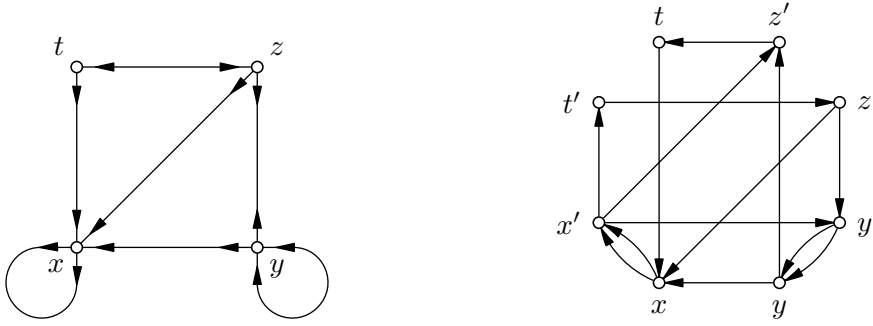

 \begin{center}
\includegraphics{pics/examples.1}%
\hspace{3cm}%
\includegraphics{pics/examples.2}%
 \end{center}
\caption{Related bidirected and skew-symmetric graphs}
\label{fig:sk-bi}
\end{figure}

Conversely, a bidirected graph $H$ with node set $V_1$, say,
determines skew-symmetric graph $G=(V,E)$ with symmetry $\sigma$ as
follows. Take a copy $\sigma(v)$ of each element $v$ of $V_1$,
forming the sets $V_2:=\{\sigma(v):v\in V_1\}$ and
$V:=V_1\sqcup V_2$. For each edge $e$ of $H$ connecting nodes
$u$ and $v$, assign two ``symmetric'' arcs $a,a'$ in $G$ so as to satisfy
(i)-(iii) above (where $u'=\sigma(u)$ and $v'=\sigma(v)$). An
example is depicted in Fig.~\ref{fig:sk-bi}.

\medskip
\noindent{\bf Remark 4.} A bidirected graph generates one
skew-symmetric graph, while a skew-symmetric graph generates a
number of bidirected ones, depending on the partition $\pi$ of $V$
that we choose in the first construction. The latter bidirected
graphs are produced from each other by the edge reversing
transformation with respect to a subset of nodes as indicated in
Remark 3 in the Introduction, so they are equivalent for us.

\medskip
A terminal set $S$ in $H$ generates the symmetric terminal set
$T:=S\sqcup\sigma(S)$ in $G$, and vice versa. One easily checks that
$(H,S)$ is inner Eulerian if and only if $(G,T)$ is such. Also there
is a correspondence between the $S$-paths in $H$ and certain
$T$-paths in $G$. More precisely, let $\tau$ be the natural mapping
of $V\cup E$ to $V_1\cup E(H)$ (where $E(H)$ is the edge set of $H$).
Each walk (cycle) $P=(v_0,a_1,v_1,\ldots,a_k,v_k)$ in $G$ induces
the sequence $\tau(P):=
(\tau(v_0),\tau(a_1),\tau(v_1),\ldots,\tau(a_k),\tau(v_k))$ of
nodes and edges in $H$.

Conversely, for a walk (cycle)
$Q=(w_0,e_1,w_1,\ldots,e_k,w_k)$ in $H$, form the sequence
$\tilde \tau(Q):=(v_0,a_1,v_1,\ldots,a_k,v_k)$ of nodes and arcs in
$G$ by the following rule:
  \begin{itemize}
\item[(R)] $v_0:=w_0$ if $e_1$ leaves $w_0$, and $v_0:=\sigma(w_0)$
if $e_1$ enters $w_0$; and for $i=1,\ldots,k$: (a) if $e_i$ leaves
$w_{i-1}$, then $a_i$ is the arc in $\tau^{-1}(e_i)$ that leaves
$w_{i-1}$, and $v_i$ is the head of $a_i$; (b) if $e_i$ enters
$w_{i-1}$, then $a_i$ is the arc in $\tau^{-1}(e_i)$ that leaves
$\sigma(w_{i-1})$, and $v_i$ is the head of $a_i$.
  \end{itemize}
(When $e_i$ is a loop, the arcs in $\tau^{-1}(e_i)$ are parallel,
and the arc $a_i$ in this set is chosen arbitrarily.)
It is not difficult to conclude that (R) provides:
  \begin{numitem}
for a walk (cycle) $Q$ in $H$,
  \begin{itemize}
 \item[(i)] $\tilde\tau(Q)$ is a walk (cycle) in $G$ and
$\tau(\tilde\tau(Q))=Q$;
 \item[(ii)] if $Q$ is edge-simple and minimal (see Remark 2 in the
Introduction), then $\tilde\tau(Q)$ is a regular path (circuit).
  \end{itemize}
  \label{eq:corresp}
  \end{numitem}
Also the walk (cycle) reverse to $Q$ determines the walk (cycle) in
$G$ symmetric to $\tilde\tau(Q)$ (up to the choice of arcs $a_i$
for loops $e_i$). The corresponding converse properties to those
in~\refeq{corresp} also take place.

Let us say that a $T$-walk $P$ from $s$ to $t$ in $G$
is {\em essential} if $t$ is different from $\sigma(s)$.
Thus, we have a natural bijection between the essential regular
$T$-paths in $G$ (considered up to parallel arc mates) and the
minimal $S$-paths in $H$. This gives
  \begin{equation} \label{eq:twice1}
\tilde\nu_{G,T}=2\nu_{H,S},
  \end{equation}
where $\tilde\nu_{G,T}$ is the maximum cardinality of a symmetric
collection of pairwise arc-disjoint essential $T$-paths in $G$.

Note also that for a terminal $s\in S$ and an $s$-set $X$ in $H$,
each edge of the cut $\delta_H(X)$ generates two arc mates in
the symmetric cut $\delta_G(X\sqcup\sigma(X))$ in $G$. Therefore,
  \begin{equation} \label{eq:twice2}
\tilde\lambda_{G,T}(s)=2\lambda_{H,S}(s)\qquad
\mbox{for each $s\in S$},
  \end{equation}
where $\tilde\lambda_{G,T}(s)$ is the minimum cardinality of a
symmetric cut in $G$ separating $\{s,s'\}$ and $T-\{s,s'\}$.

In view of relations~\refeq{twice1} and~\refeq{twice2},
Theorem~\ref{tm:main} is equivalent to the following

  \begin{corollary} \label{cor:skew}
For a skew-symmetric graph $G=(V,E)$ with a symmetric set
$T$ of terminals, if $(G,T)$ is inner Eulerian, then
$\tilde\nu_{G,T}=\frac12\sum_{s\in S}\tilde\lambda_{G,T}(s)$.
   \end{corollary}

In the capacitated case, we are given a symmetric function
$c:E\to\Zset_+$ of arc capacities in a skew-symmetric graph
$G=(V,E)$ with a symmetric set $T=S\sqcup S'$ of terminals.
By an integer symmetric free multiflow (or, briefly, an
IS-{\em multiflow}) in the network $(G,T,c)$ we mean a collection
$F$ of integer flows $f_{st}$
for the ordered pairs $(s,t)$ of distinct terminals in $S$ such
that: (a) $f=f_{st}$ is a flow from $\{s,s'\}$ to $\{t,t'\}$
(i.e., $\diver_f(v)$ is nonnegative for $v=s,s'$, nonpositive for
$v=t,t'$, and 0 otherwise); (b) each $f_{st}$ is symmetric to
$f_{ts}$; and (c) $F$ is $c$-{\em admissible}, i.e.,
   $$
\sum\nolimits_{st} f_{st}(e)\le c(e) \qquad
  \mbox{for each $e\in E$}.
  $$
The (total) {\em value} $\valu(F)$ of $F$ is $\sum_{st}\valu(f_{st})$.
The problem is:
  \begin{itemize}
\item[(PS)]
{\em Given an inner Eulerian network $(G,T,c)$, where $G$ is a
skew-symmetric graph and $T$ and $c$ are symmetric, find a {\em
maximum IS-multiflow}, i.e., an IS-multiflow $F$ maximizing
$\valu(F)$.}
  \end{itemize}

To see how this problem is related to~(P), let
$(G,T)$ correspond to $(H,S)$, where $H$ is bidirected. Let
$\hat c$ be the corresponding capacity function in $H$,
i.e., $\hat c(e)=c(a)$ for an edge $e\in E(H)$ and its images $a,a'$
in $G$. The inner Eulerianness of $(G,T,c)$ implies that of
$(H,S,\hat c)$, and vice versa. Given an IS-multiflow $F$ in
$(G,T,c)$, represent each flow $f_{st}$ in the path packing form:
  \begin{numitem}
$f_{st}=\alpha_1\chi^{P_1}+\ldots+\alpha_k\chi^{P_k}$, where
$\alpha_i=\alpha(P_i)\in\Zset_+$ and $P_i$ is a circuit or a
(simple) path from $\{s,s'\}$ to $\{t,t'\}$.
  \label{eq:fst}
  \end{numitem}
We assume that the representation of each flow $f_{st}$ is symmetric
to that of $f_{ts}$. Then the
set $\Pscr$ of (essential) $T$-paths in these representations is
symmetric, with $\alpha(P)=\alpha(P')$ for each $P\in\Pscr$, and
we have $\valu(F)=\sum(\alpha(P):P\in\Pscr)$. Now each pair
$P,P'\in\Pscr$ of path mates determines an $S$-walk $\hat P$
in $H$ (considered up to reversing), and taking together these paths
$\hat P$ with weights $\alpha(P)$, we obtain a multiflow
$\Fscr$ in $(H,S,\hat c)$ satisfying $\valu(\Fscr)=\frac12\valu(F)$.

Conversely, let $\Fscr=(\hat\Pscr,\hat\alpha)$ be an integer
multiflow in $(H,S,\hat c)$, where $\hat\Pscr$ consists of
$S$-walks. One may assume that for each edge $e$ of $H$, no path
$\hat P\in\hat\Pscr$ traverses $e$ twice in the same direction
(for otherwise one can remove a cycle from $\hat P$). Then each
$s$--$t$ walk $\hat P$ determines an arc-simple directed walk $P$
from $\{s,s'\}$ to $\{t,t'\}$ and its mate $P'$ from $\{t,t'\}$ to
$\{s,s'\}$ in $G$.
Assign $\alpha(P):=\alpha(P'):=\hat\alpha(\hat P)$. Let $f_{st}$ be
the sum of functions $\alpha\chi^P$ over the obtained walks $P$ from
$\{s,s'\}$ to $\{t,t'\}$. Then $f_{ts}$ is symmetric to $f_{st}$.
These flows form an IS-multiflow $F$ in
$(G,T,c)$ satisfying $\valu(F)=2\valu(\Fscr)$.

Thus, problems~(PS) and~(P) (regarding $H,S,\hat c$) are reduced to
each other. In the next section we devise an efficient algorithm
for finding an optimal solution to~(PS) and then explain that it
can be transformed into an optimal solution to the corresponding
instance of~(P) without increasing the time bound.

\section{{\Large\rm Algorithm}} \label{sec:alg}

In this section we describe an algorithm to solve
problem~(PS) and estimate its complexity. We use
terminology and facts from the previous section.

Let $(G=(V,E),T,c)$ be an inner Eulerian skew-symmetric network. As
before, we represent the terminal set $T$ as $S\sqcup S'$ and
associate with $(G,T,c)$ the corresponding bidirected network
$(H,S,\hat c)$. One may assume that no arc in $G$ connects a
pair of terminal mates. Also if $G$ has an arc entering a terminal
$s\in S$, then replacing its head $s$ by $s'$ and symmetrically
replacing the tail $s'$ of the symmetric arc $e'$ by $s$ does not
affect the problem in essence. So we may
assume that $\degin(s)=0$ for each terminal
$s\in S$ in $G$. Then any flow from $\{s,s'\}$ to $\{t,t'\}$,
where $s,t\in S$, is essentially a flow from $s$ to $t'$, and its
symmetric flow is a flow from $t$ to $s'$; this property will simplify
technical details in our construction. In terms of $H$, the latter
assumption says that each edge incident with a terminal $s$ in $H$
leaves $s$ (cf. Remark~3 in Section~\ref{sec:intr}).

The algorithm uses a recursion analogous to that
in~\cite{kar-79}, and the case $|S|=3$ is the base in it.
We first consider this special case (which generalizes the case
$|S|=2$).

  \subsection{{\large\rm Case $|S|=3$.}}  \label{ssec:three}
The algorithm for this case uses one auxiliary skew-symmetric
graph $G_1=(V,E_1)$. It is obtained from $G$ by
adding, for each pair $v,v'$ of inner node mates, four {\em
auxiliary arcs} connecting $v$ and $v'$: two arc mates
$a_v,a'_{v}$ going from $v$ to $v'$ and two arc mates
$a_{v'},a'_{v'}$ from $v'$ to $v$, regardless of the existence of
such arcs in $G$. This $G_1$ corresponds to the bidirected graph
$H_1$ obtained from $H$ by adding two auxiliary loops at each inner
node $v$, one leaving $v$ (twice) and the other entering $v$.

The algorithm consists of three stages. Let
$S=\{s_i,s_2,s_3\}$.

\medskip
At {\em Stage 1}, we apply the algorithm for inner Eulerian graphs
from~\cite{IKN-98} to find a maximum integer free multiflow in the
underlying undirected graph $\bar H$ for $H$ having the same set
$S$ of terminals and the same capacities $\hat c$. It runs in
$O(\phi(V,E))$ time (since $|S|=O(1)$) and outputs (simple)
$S$-paths $\bar P_1,\ldots,\bar P_k$ in $\bar H$ and weights
$\alpha_1,\ldots,\alpha_k\in\Zset_+$ satisfying the packing
condition w.r.t. $\hat c$. (Recall that $\phi$ is a bound for the
applied max flow algorithm; we assume $\phi=\Omega(VE)$).
It also outputs pairwise disjoint
$s_i$-sets $\bar X_i$, $i=1,2,3$, such that for each $i$, the sum
of weights $\alpha_j$ of paths $\bar P_j$ connecting $s_i$ and
$S-\{s_i\}$ is equal to $\hat c(\delta_H(\bar X_i))$.
However, some pairs of consecutive edges in $\bar P_j$ may be
non-transit in $H$, i.e., $\bar P_j$ is not necessarily a path in
$H$.

\medskip
At {\em Stage 2}, we transform $\bar P_1,\ldots,\bar P_k$ into
paths in the auxiliary bidirected graph $H_1$. More precisely,
for each $\bar P_i=(v_0,e_1,v_1,\ldots,e_q,v_q)$ and for each
non-transit pair $e_j,e_{j+1}$ in it, if both edges
$e_j,e_{j+1}$ enter (leave) $v_j$, then the element $v_j$ of
$\bar P_i$ is replaced by the string $v_j,\ell,v_j$, where $\ell$
is the auxiliary loop leaving (resp. entering) $v_j$. This results
in minimal $S$-paths $\tilde P_1,\ldots,\tilde P_k$ in $H_1$.

\begin{figure}[tb]
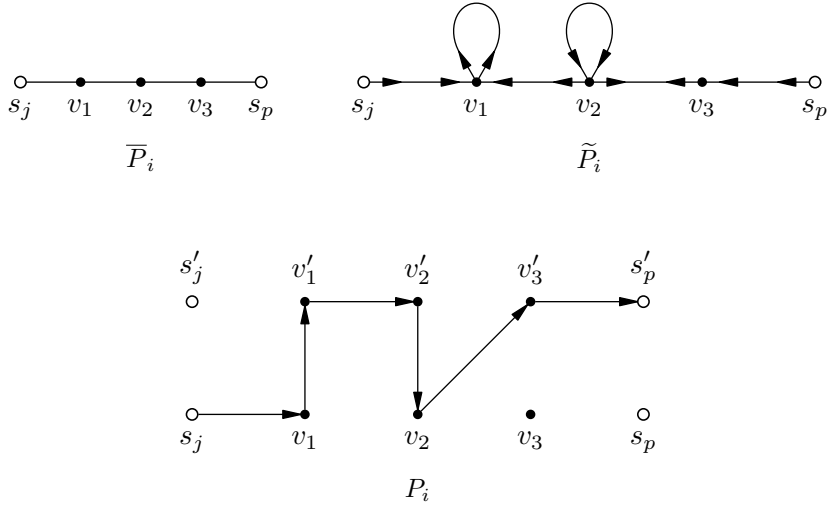

\begin{center}
\includegraphics{pics/paths.1}%
\hspace{1cm}%
\includegraphics{pics/paths.2}
\vskip 1cm
\includegraphics{pics/paths.3}
\end{center}
\caption{A path $\bar P_i$ and its images in the graphs $H_1$ and
$G_1$}
\label{fig:Pi}
\end{figure}

\begin{figure}[tb]
\begin{center}
\includegraphics{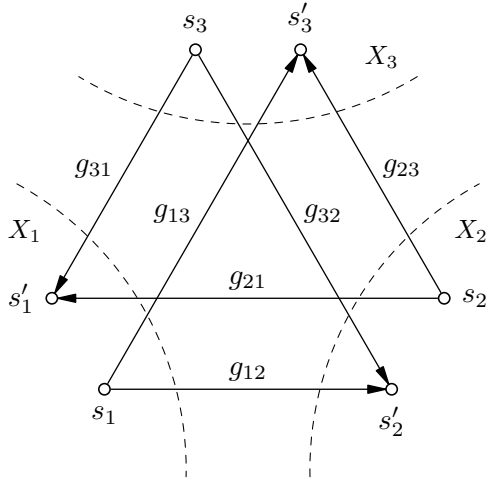}%
\end{center}
\caption{Flows $g_{ij}$ in $G_1$}
\label{fig:3term}
\end{figure}

Each path $\tilde P_i$ and its reverse one are then lifted to $G_1$
(by the method explained in Section~\ref{sec:skew}), giving regular
$T$-paths $P_i,P'_i$ symmetric to each other.
(Figure~\ref{fig:Pi} illustrates paths $\bar P_i,\tilde P_i,P_i$.)
For each pair $s_j,s_p$ ($j\ne p$),
the functions $\alpha_i\chi^{P_i}$ or $\alpha_i\chi^{P'_i}$ for the
paths from $s_j$ to $s'_p$ are added up, forming $s_j$--$s'_p$ flow
$g_{jp}$. This gives a symmetric collection of six integer flows in
$G_1$; see Figure~\ref{fig:3term}. The $\hat c$-admissibility of
the above multiflow in $\bar H$ and the fact that each path $P_i$ is
regular imply that the total flow though each arc
$e$ of $G$ does not exceed $c(e)$. Also the fact that the cuts
$\delta_H(\bar X_i)$ are saturated implies that
  \begin{numitem}
for $i=1,2,3$, the arcs in $\deltaout(X_i)$ are saturated by
$g_{i,i-1}+g_{i,i+1}$, and symmetrically, the arcs in
$\deltain(X_i)$ are saturated by $g_{i-1,i}+g_{i+1,i}$,
  \label{eq:satur}
  \end{numitem}
where $X_i:=\bar X_i\sqcup \bar X'_i$ and the indices are taken
modulo 3. So the IS-multiflow consisting of these six flows has
maximum value.

\medskip
At {\em Stage 3}, we improve the above flows $g_{ij}$ in
$G_1$ by reducing their values on the auxiliary arcs to zero,
eventually obtaining the desired multiflow in $G$. In view
of~\refeq{satur}, for $i=1,2,3$, one may assume that $g_{i-1,i+1}$ and
$g_{i+1,i-1}$ take zero values on all arcs of the subgraph $\langle
X_i\rangle$ of $G_1$ induced by $X_i$.

Take the residual capacities $\Delta(e):=c(e)-
\sum_{ij}g_{ij}(e)$ of arcs $e\in E$. The divergency of $c$
(w.r.t. $E$) and of each $g_{ij}$ (w.r.t. $E_1$) at any inner node is
zero, therefore,
  \begin{equation} \label{eq:divDelta}
\diver_\Delta(v)=\sum\nolimits_{ij} g_{ij}(v,v') \qquad
  \mbox{for each $v\in V-T$},
  \end{equation}
where $g(v,v')$ denotes $g(a_v)+g(a'_v)-g(a_{v'})-g(a'_{v'})$
(recall that $a_v,a'_v$ are the auxiliary arcs from $v$ to $v'$).
The function $\Delta$ on $E$ is nonnegative, integer-valued and
symmetric. Also~\refeq{divDelta} and
  \begin{equation} \label{eq:bal}
 g_{ij}(v,v')=-g_{ij}(v',v)=g_{ji}(v,v')\qquad
         \mbox{for each $v\in V-T$}
  \end{equation}
imply that $\diver_\Delta(v)$ is even for each $v\in V-T$.
Hence we can extend $\Delta$ to the auxiliary arcs so as to obtain
an IS-flow in $(G_1,T)$. (The extended $\Delta$ satisfies
$\Delta(v,v')+\sum_{ij}g_{ij}(v,v')=0$ for each $v\in V-T$.)

Notice that $\Delta(e)=0$ for each arc $e$ in the cut
$\delta(X_i)$, $i=1,2,3$, by~\refeq{satur}. Therefore, the
restriction $\Delta_i$ of $\Delta$ to the set $A_i$ of arcs of
the subgraph $\langle X_i\rangle$ is an IS-flow from $s_i$ to
$s'_i$. In its turn, the
restriction $\Delta_0$ of $\Delta$ to the set $A_0$ of arcs with
both ends in $W:=V-(X_1\cup X_2 \cup X_3)$ is an
integer symmetric circulation in the subgraph
$(W,A_0)$. (Recall that the sets $X_1,X_2,X_3$ are pairwise
disjoint.)

We start with getting rid of nonzero arc values of the above flows
on the auxiliary arcs within the subgraph $\langle X_1\rangle$.
To this aim, apply the integer symmetric decomposition procedure
to $\Delta_1$ (cf.~\refeq{sdecom})
to represent it as the sum of integer
$s_1$--$s'_1$-flows $h,h'$, where $h'$ is symmetric to $h$.
Combine $g:=g_{12}+g_{13}+h$ and $g':=g_{21}+g_{31}+h'$
(where $h,h'$ are formally extended by zeros on $E_1-A_1$).
Then $g$ is an integer flow from $s_1$ to $S'$, and $g'$ is the flow
from $S$ to $s'_1$ symmetric to $g$. Also
  $$
  g(v,v')=g_{12}(v,v')+g_{13}(v,v')+h(v,v')=0\qquad
       \mbox{for each $v\in X_1-\{s_1,s'_1\}$},
  $$
in view of $h(v,v')=\frac12\Delta_1(v,v')$, $\Delta(v,v')+
\sum_{ij}g_{ij}(v,v')=0$ and $g_{23}(e)=g_{32}(e)=0$ for all
$e\in A_1$. So we can
reduce $g,g'$ to zero on all auxiliary arcs in $\langle
X_1\rangle$. Now using standard flow decomposition, we represent
the new flow $g$ as the sum of three integer flows $f_1,f_2,f_3$,
from $s_1$ to $s'_1$, from $s_1$ to $s_2'$, and from $s_1$ to
$s_3'$, respectively. Note that $g(e)=0$ for each
$e\in\deltain(X_1)$ implies that $f_1$ is zero on all arcs of the cut
$\delta(X_1)$. Update $g_{12}:=f_2$ and $g_{13}:=f_3$; the flows
$g_{21}$ and $g_{31}$ are updated symmetrically. Then the resulting
four flows together with the remaining flows $g_{23},g_{32}$
satisfy~\refeq{satur} as before (thus forming a maximum
IS-multiflow) and take zero values on the auxiliary arcs in
$\langle X_1\rangle$, as required. Do similarly for $X_2$ and
$X_3$.

The task of improving the flows within the subgraph
$\langle W\rangle=(W,A_0)$ is a bit more involved.
First of all we modify $g_{12}$ (and $g_{21}$) so as to get
  \begin{equation} \label{eq:zero}
 g_{12}(v,v')+g_{23}(v,v')+g_{31}(v,v')=0\qquad
            \mbox{for each $v\in W$}
  \end{equation}
(this situation is technically simpler).
This is done by decomposing the above-mentioned symmetric
circulation $\Delta_0$ in $\langle W\rangle$ into the sum of
an integer circulation $\omega$ and its symmetric circulation
$\omega'$ and then by updating $g_{12}:=g_{12}+\omega$ and
$g_{21}:=g_{21}+\omega'$ (with $\omega,\omega'$
extended by zeros to $E_1-A_0$). Then the equality $\Delta(v,v')+
\sum_{ij}g_{ij}(v,v')=0$ provides~\refeq{zero}.

The process of improving the flows within $\langle W\rangle$
consists of $O(W)$ iterations (the idea is borrowed from the
algorithm for digraphs in~\cite{IKN-98}). At a current iteration,
we choose a node $v\in W$ where some $g_{ij}(v,v')$ is nonzero.
W.l.o.g, one may assume that $g_{12}(v,v')>0$ and $g_{13}(v,v'),
g_{23}(v,v')\le 0$. Let $r_0:=|g_{13}(v,v')|$ and $r_1:=
|g_{23}(v,v')|$; then $g_{12}(v,v')=r_0+r_1$, by~\refeq{zero}. Let
$B$ be the set of (four) auxiliary arcs connecting $v$ and $v'$.

First of all we represent $g_{12}$ as the sum of two integer
$s_1$--$s'_2$ flows $g_0,g_1$ such that $g_0(v,v')=r_0$ and
$g_1(v,v')=r_1$. To do so, replace $B$ by new terminals
$t,t_0,t_1$ and arcs $a=(t,v')$, $a_0=(v,t_0)$ and $a_1=(v,t_1)$,
and add an arc $b$ from $s'_2$ to $s_1$. Define $h(a_i):=r_i$, $i=0,1$,
$h(a):=r_0+r_1$, $h(b):=\valu(g_{12})$ and $h(e):=g_{12}(e)$ for
the remaining arcs $e$. This turns $h$ into a flow from $t$ to
$\{t_0,t_1\}$, and we decompose it into the sum of integer flows
$h_0,h_1$, from $t$ to $t_0$ and from $t$ to $t_1$, respectively.
These $h_0,h_1$ determine the desired $g_0,g_1$ in a natural way.

Combine $f:=g_0+g_{13}$. Then $f(v,v')=g_0(v,v')+g_{13}(v,v')=0$.
Update $f(e):=0$ for each $e\in B$ and decompose the updated
$s_1$--$\{s'_2,s'_3\}$ flow $f$ into the sum of integer flows
$f_0,f_{13}$, from $s_1$ to $s'_2$ and from $s_1$ to
$s'_3$, respectively. Then $\valu(f_0)=\valu(g_0)$ and
$\valu(f_{13})=\valu(g_{13})$.

Doing similarly for the flow $g_{23}$ and the flow $g'_1$ symmetric
to $g_1$ (which have the source $s_2$ in common), we obtain
corresponding $s_2$--$s'_3$ flow $f_{23}$ and $s_2$--$s'_1$ flow
$f'_1$. Finally, update $g_{12}:=f_0+f_1$ (where $f_1$ is symmetric
to $f'_1$), $g_{13}:=f_{13}$ and $g_{23}:=f_{23}$. The updated flows
$g_{ij}$ together with their symmetric ones satisfy $g_{ij}(e)=0$
for each arc $e\in B$.

Then we choose a next pair of node mates in $W$, and so on. Upon
termination of the process, the resulting flows $g_{ij}$ take zero
values on all auxiliary arcs, and it is easily seen from the
construction that $\valu(g_{ij})$ preserves for all pairs $ij$. So their
restrictions to $E$ form a maximum IS-multiflow in $(G,T,c)$, as
required.

\medskip
The above algorithm runs in $O(\phi(V,E))$ time plus the time
needed to perform $O(W)$, or $O(V)$, flow decompositions
during the iterative process at Stage 3 (the other
operations including those in $O(1)$ symmetric decompositions take
$O(VE)$ time). Each of these decompositions is applied to a flow
with $O(1)$ sources and sinks, and we use the procedure
in Section~\ref{sec:decom} to implement it in $O(E\log(2+V^2/E))$
time. This gives the bound $O(VE\log(2+V^2/E))$ for the
six (or four) terminal case.

  \subsection{{\large\rm General case.}}  \label{ssec:gen}
We now describe the algorithm for an arbitrary $|S|\ge 4$. It is
based on a recursive network partition approach.

For a {\em current} inner Eulerian skew-symmetric network $N=(G,T,c)$,
with $T=S\sqcup S'$, the {\em network partition procedure} partitions
$S$ into two sets $S_1,S_2$ such that $|S_1|=\lceil |S|/2\rceil$ and
$|S_2|=\lfloor |S|/2\rfloor$
and finds a symmetric subset $X\subset V$ with $X\cap T=
S_1\sqcup S'_1$ whose induced cut $\delta(X)$ has minimum capacity
$c(\delta(X))$. This is done by finding a minimum capacity cut
$\delta(Y)$ with $Y\cap T=S_1\sqcup S'_1$ in the underlying
undirected network for $(G,c)$, and by making the symmetrization
$X:=Y\cup Y'$ (relying on $c(\delta(Y\cup Y'))+c(\delta(Y\cap Y'))
\le c(\delta(Y))+c(\delta(Y'))=2c(\delta(Y))$).

Next we shrink the subgraph $\langle V-X\rangle$ of $G$ into two new
({\em extra}) terminals $t_1,t'_1$,
making each arc in $\deltaout(X)$ enter $t'_1$, and each arc in
$\deltain(X)$ leave $t_1$.
Similarly, $\langle X\rangle$ is shrunk into extra terminals
$t_2,t'_2$, each arc in $\deltain(X)$ becomes entering $t'_2$ and
each arc in $\deltaout(X)$ becomes leaving $t_2$. This produces two
smaller inner Eulerian networks $N_i=(G_i=(V_i,E_i),T_i,c_i)$ with
$T_i=S_i\cup S'_i\cup\{t_i,t'_i\}$, $i=1,2$, satisfying
  \begin{equation} \label{eq:TVE}
 |T_i|\le 4|T|/5,\quad |V_i|\le |V|, \quad |E_i|\le |E|,
  \quad\mbox{and} \quad |V_1|+|V_2|=|V|+4
  \end{equation}
(since $|T_1|=8$ when $|T|=10$). Also for $X_i:=\{t_i,t'_i\}$, the cut
$\delta(X_i)$ of $G_i$ has minimum capacity among the cuts
separating $\{t_i,t'_i\}$ and $T_i-\{t_i,t'_i\}$.

One application of the network partition procedure, to a current
$N$, takes one minimum cut computation, so it runs in
$O(\phi(V,E))$ time.

Let $F_i$ be a (recursively found) maximum free IS-multiflow in
$N_i$. The {\em aggregation procedure} transforms $F_1,F_2$ into a
maximum free IS-multiflow $F$ in $N$. The
flows in $F_i$ going from $S_i$ to $t'_i$ are combined into one
(multisource) flow $f_i$ from $S_i$ to $t'_i$, and symmetrically,
the flows from $t_i$ to $S'_i$ are combined into one flow $f'_i$.
By the maximality of $F_i$ and the minimality of $c_i(\delta(X_i))$,
$f_i$ saturates $\deltain(X_i)$ and $f'_i$ saturates
$\deltaout(X_i)$. We glue together (the images of) $f_1$ and
$f'_2$, obtaining $S_1$--$S'_2$ flow $f$ in $N$, and do
symmetrically for $f_2,f'_1$, obtaining $f'$. These $f,f'$ are
decomposed symmetrically into a symmetric collection of
integer one-source-one-sink flows.
Then the flows formed from $f,f'$ together with the remaining flows in
$F_1,F_2$ connecting pairs of terminals in $T_1$ or in $T_2$ give
the desired $F$. (The maximality of $F$ follows from the fact that
for each $s\in S_1$, the total value of flows in $F_1$ leaving $s$
or entering $s'$ is equal to the minimum capacity of a cut in $N_1$
separating $\{s,s'\}$ and $T_i-\{s,s'\}$, and similarly for $S_2$.
The above construction maintains such an equality for $F$
and each $s\in S$.)

At the bottom level ($|S|=3$), we apply the algorithm described
in~\ref{ssec:three}.

One application of the aggregation procedure, to current $N_1,N_2$,
takes $O(S_1E_1+S_2E_2)$ time to create the flows $f_1,
f_2$ as above plus $O(VE)$ time to decompose $f$, or $O(VE)$ time in
total (in view of~\refeq{TVE}).

\medskip
It remains to explain that the resulting multiflow $F$ in the
initial network $N$ can be efficiently transformed into a maximum
integer free multiflow in the corresponding bidirected network
$(H,S,\hat c)$. We show that $O(VE\log T)$ time is sufficient to
create from $F$ a corresponding symmetric collection
$(\Pscr,\alpha)$ of weighted $T$-paths in $N$; these paths
determine weighted $T$-walks in $H$ forming an optimal solution to
problem~(P) with $(H,S,\hat c)$, by the relationship explained in
Section~\ref{sec:skew}. We assume that each flow $f$ in $F$ is
explicitly given only within its support $\supp(f):=\{e\in E:
f(e)\ne 0\}$.

Let $\Tscr$ be the binary rooted tree formed by all networks
arising during the recursion, with the
natural ordering on them. The height of $\Tscr$ (or the depth of
the recursion) is $O(\log T)$, in view of the first inequality
in~\refeq{TVE}. For a network $\tilde N$ in $\Tscr$, let $A(\tilde
N)$ be the set of terminals from the initial $T$ that are contained in
$\tilde N$, and $F(\tilde N)$ the set of flows in $F$ with both
terminals in $A(\tilde N)$. We use the fact that for incomparable
$N',N''$ in $\Tscr$, the supports of flows in $F(N')$ are disjoint
from those in $F(N'')$. (Indeed, for the closest common
predecessor $\tilde N$ of $N',N''$, the minimum cut found by the
network partition procedure for $\tilde N$ separates $A(N')$ and
$A(N'')$ and is saturated by the flows not in $F(N')\cup F(N'')$.)

We proceed as follows. For each non-leaf network $\tilde N$ with
children $N_1,N_2$, combine the flows in $F$ with the source in
$A(N_1)$ and the sink in $A(N_2)$ into one multiterminal flow
$\tilde f$ (in the initial network), and then decompose $\tilde f$
into a set $\Pscr(\tilde N)$ of weighted paths from $A(N_1)$ to
$A(N_2)$ (the circuits appeared in the decomposition are removed).
This takes $O(V\supp(\tilde f))$ time.
Taken together, the sets $\Pscr(\tilde N)$,
their symmetric sets and corresponding paths appeared by
decomposing the flows in $F$ having both terminals in one leaf
network, constitute the desired symmetric collection
$(\Pscr,\alpha)$. To estimate the complexity, consider the networks
$\tilde N$ at height $i$ in $\Tscr$. They are incomparable, so
the supports of flows $\tilde f$ as above in them are pairwise
disjoint. Hence to form the sets $\Pscr(\tilde N)$ for these
$\tilde N$ takes $O(VE)$ time in total. This gives the bound
$O(VE\log T)$ for the whole procedure, as declared.

  \subsection{{\large\rm Complexity of the algorithm}}
            \label{ssec:time}

We show that the above algorithm runs in $O(VE\,\ell(V,E)\log T)$
time, where $\ell(V,E):=\max\{1,\ln(V^2/E)\}$, assuming
$\phi(n,m)=O(nm\log(2+n^2/m))$ (as in Goldberg--Tarjan's max flow
algorithm). We use induction on the height $h(T)$ of the binary
tree $\Tscr$ (it depends only on $|T|$). When $h(T)=0$
(i.e. $|T|=6$), the required time bound was shown
in~\ref{ssec:three}.

Let $h(T)\ge 1$ and let $N_1,N_2$ be the children of $N$ in
$\Tscr$. For $i=1,2$, we have $h(T_i)\le h(T)-1$, and by induction
the time $\tau_i$ of the algorithm to solve the problem for $N_i$
is bounded from above as
  $$
 \tau_i\le Cn_im_i\ell(n_i,m_i)\log T_i
  $$
for some appropriately chosen constant $C>0$ (specified later).
Here $n_i:=|V_i|$ and $m_i:=E_i$, keeping notation
from~\ref{ssec:gen}. The network partition and aggregation
procedures applied to $N$ take time $O(\phi(n,m))$ and $O(nm)$,
respectively, or $Dnm\ell(n,m)$ time together, where $D$ is some
constant $>0$, $n:=|V|$ and $m:=|E|$. Therefore, the time $\tau$ to
solve the problem for $N$ is estimated as
  \begin{equation}  \label{eq:time1}
\tau\le C(n_1m_1\ell(n_1,m_1)\log T_1+n_2m_2\ell(n_2,m_2)\log T_2)
  +Dnm\ell(n,m).
  \end{equation}

We have $\ell(n_i,m_i)\le\ell(n,m_i)$ (since $n_i\le n$) and
$m_i\ell(n,m_i)\le m\ell(n,m)$ (this follows from $m_i\le m$
and from $\ln a>\frac{1}{b}\ln(ab)$ for $b>1$ and $\ln a>1$). Also
$n_1+n_2=n+4$ and $\log T_i\le\log T-\log\frac54$,
by~\refeq{TVE}. Then~\refeq{time1} implies
  \begin{equation} \label{eq:time2}
\tau\le C\left( nm\log T-nm\log\frac54+4m\log T-4m\log\frac54\right)
           \ell(n,m)+Dnm\ell(n,m).
  \end{equation}

Since $n\ge|T|$ and $\frac12 nm\log \frac54$ grows faster
than $4m\log n$, one can choose constants $n_0$ and $C$
(depending on $D$) such that the right hand side value
in~\refeq{time2} becomes smaller than $Cnm\ell(n,m)\log T$ for any
$n>n_0$. (For the networks with $|V|\le n_0$, the problem is
solved in $O(E)$ time.) This yields the desired time bound.
  \begin{theorem} \label{tm:time}
A maximum IS-multiflow (resp. a maximum integer free multiflow) in
an inner Eulerian skew-symmetric (resp. bidirected) network
$(G=(V,E),T,c)$ can be found in $O(VE\log T\log(2+V^2/E))$ time.
  \end{theorem}

\section{{\Large\rm Fast Flow Decomposition}} \label{sec:decom}

For a fixed $k\in\Zset_+$, we consider the problem:
  \begin{itemize}
\item[(D)] {\em Given a flow (integer flow) $f$ from $S$ to $T$,
with $|S|+|T|=k$, in a digraph $G=(V,E)$, find a decomposition
$f=\sum(f_{st}: s\in S,t\in T)$, where each $f_{st}$ is a flow (resp.
integer flow) from $s$ to $t$,}
  \end{itemize}
and show the following (allowing parallel arcs in $G$ and assuming
$|V|=O(E)$).
  \begin{theorem}  \label{tm:decom}
(D) can be solved in $O(E\log(2+V^2/E))$ time.
  \end{theorem}

Note that when $G$ is acyclic, a decomposition
(into one-source-one-sink flows or into weighted paths) of any flow
in $G$ is carried out
in $O(E)$ time by using a topological sorting of the nodes.
Sleator and Tarjan~\cite{ST-83} showed that any flow $f$ in an
arbitrary digraph can be decomposed, in $O(E\log V)$ time, into a
circulation and a flow whose support induces an acyclic
subgraph of $G$ (so a decomposition of $f$ into one-source-one-sink
flows can be found with the same complexity $O(E\log V)$). The
algorithm in~\cite{ST-83} uses sophisticated computational tools,
so-called dynamic trees.

Our approach to solve~(D) is based on a node splitting technique
and uses only simple data structures. Let $\Pi$ be the set of pairs
$st$ with $s\in S$ and $t\in T$.

In the beginning of the algorithm, we delete from $G$ the arcs $e$
with $f(e)=0$. Also we sort the nodes $v$ by increasing their degrees
$\deg(v)$. (This takes $O(E)$ time.) The algorithm applies $|V|$
iterations.

At each iteration, we choose a node $v$ with $\deg(v)$ minimum in the
current graph $G=(V,E)$. First of all we scan the arcs incident with
$v$ to select parallel arcs among them. Each tuple of parallel arcs is
merged into one arc (and the flows on these are added up). The node
degrees and the ordering on $V$ are updated accordingly. (This
preliminary stage is performed in $O(\deg(v))$ time. As a result, the
degree of $v$ becomes less than $2|V|$.) Then we make at most $\deg(v)$
splittings at $v$.

More precisely, at a current step
of the iteration, we choose an arc entering $v$ and an arc leaving
$v$, say, $e=(u,v)$ and $e'=(v,w)$. If $e$ or $e'$ is a loop, we
simply delete it from $G$. Otherwise define $\eps:=\min\{f(e),f(e')
\}$. The {\em splitting-off operation} applied to $(e,v,e')$
creates a new arc $e''$ from $u$ to $w$, assigns $f(e''):=\eps$,
updates $f(e):=f(e)-\eps$ and $f(e'):=f(e')-\eps$, and deletes from
$G$ the arc (or arcs) for which the new value becomes zero;
it takes $O(1)$ time. The ordering on $V$ is updated accordingly (in
$O(1)$ time). Clearly the operation maintains both the divergency
at each node and the flow integrality (when the original flow is
integer). Also $\diver(v)$ decreases and the number of all arcs does
not increase.

At the next step of the iteration, the operation is applied to
another pair of arcs, one entering and the other leaving $v$, and
so on until such pairs no longer exist. After that, if $v\not\in
S\cup T$, then $v$ is removed from $G$ (as $\diver_f(v)=0$ implies
$\deg(v)=0$).

At the next iteration, we again choose a vertex where the
current degree is minimum, and so on. One can see that after $|V|$
iterations, each arc of the resulting graph $\bar G$ goes from a
source $s\in S$ to a sink $t\in T$. The decomposition $\bar D=
\{f_{st}:st\in\Pi\}$ for the resulting $f$ in $\bar G$ is
trivial: $f_{st}(s,t):=f(s,t)$ and $f_{st}(e):=0$ for $e\ne (s,t)$
(letting $f_{st}:\equiv 0$ if the arc $(s,t)$ does not exist in
$\bar G$).

Now going in the reverse order and applying the corresponding {\em
restoration procedure} reverse to the splitting-off one, we
transform $\bar D$ into the desired decomposition of the initial
flow. More precisely, consider a current graph $G$ and the arcs
$e=(u,v)$, $e'=(v,w)$, $e''=(u,w)$ as above, and let $f_{st}$,
$st\in\Pi$, be the flows already obtained for the graph $G'$
formed from $G$ by the splitting-off operation w.r.t. $(e,v,e')$.
For each $st\in\Pi$, add $f_{st}(e'')$ to $f_{st}(e)$ and to
$f_{st}(e')$ and then delete $e''$.
(The backward iteration concerning $v$ finishes with restoring the
corresponding tuples of parallel arcs incident with $v$ and assigning,
in a due way, the flows $f_{st}$ on these arcs.)
Eventually, we obtain the desired decomposition $\{f_{st}:
st\in\Pi\}$ of the initial $f$. (Strictly speaking, we have
$g:=\sum_{st}f_{st}\le f$ and $\diver_{f-g}(v)=0$ for all $v\in V$;
so one should add the circulation $f-g$ to one of the flows
$f_{st}$.)

\medskip
Next we estimate complexity of the above algorithm. Let $v_1,v_2,
\ldots,v_{|V|}$ be the sequence of nodes in the splitting-off
process. Since $|\Pi|=O(1)$, the restoration process is only $O(1)$
times slower than the splitting-off one. (This is just where we
essentially use the condition that $f$ has $O(1)$ terminals.)
Using this fact, one can conclude that the algorithm runs in
$O(E+\Delta)$ time for the initial $E$, where
$\Delta:=\deg^\ast(v_1)+\ldots+\deg^\ast(v_{|V|})$
and $\deg^\ast(v)$ denotes the degree of $v$ at the beginning of
splitting at $v$. Each iteration $i$ in the former process does not
increase the number of arcs of the current graph and decreases the
number of nodes by one, unless $v_i\in S\cup T$. So
$\deg^\ast(v_{i+1})$ is at most $2|E|/(|V|-i)$. Summing up the
latter numbers over $i$, we obtain $\Delta=O(E\log V)$, which is
worse than the time bound in Theorem~\ref{tm:decom}.

However, we can estimate $\Delta$ more carefully, by using the
inequality $\deg^\ast(v_i)<2(|V|-i+k+1)$ (provided by merging
parallel arcs incident with $v_i$).
For any integer $1\le\lambda\le |V|$, apply the first
bound on $\deg^\ast(v_i)$ for $i=1,\ldots,|V|-\lambda$, and the
second bound for $i=|V|-\lambda+1,\ldots,|V|$. This gives
   $$
\Delta\le 2|E|\left(\frac{1}{|V|}+\frac{1}{|V|-1}+\ldots+
  \frac{1}{\lambda+1}\right)+2\lambda(\lambda+k),
   $$
or $\Delta=O(E\log(V/\lambda)+\lambda^2)$. Now taking $\lambda:=
\min\{|V|,\lceil\sqrt{|E|}\rceil\}$, we obtain $\Delta=
O(E\log(2+V^2/E))$, and the theorem follows.

\section{{\Large\rm Fast Skew-Symmetric Flow Decomposition}}
                          \label{sec:sdecom}

In this section Theorem~\ref{tm:decom} is extended to
(skew-)symmetric flows. For a fixed $k\in\Zset_+$, we consider the
problem:
  \begin{itemize}
 \item[(DS)] {\em Given an integer symmetric flow $f$ from
$S=\{s_1,\ldots,s_k\}$ to $S'=\sigma(S)$ in a skew-symmetric graph
$G=(V,E)$, find a decomposition of $f$ of the form
  \begin{equation} \label{eq:decs}
f=\sum\nolimits_{1\le i\le j\le k}(f_{ij}+f'_{ij}),
  \end{equation}
where each $f_{ij}$ is an integer flow from $s_i$ to $s'_j$ and
$f'_{ij}$ is symmetric to $f_{ij}$.}
  \end{itemize}

Note that $f'_{ij}$ is a flow from $s_j$ to $s'_{i}$. So
in the above decomposition, for $i,j\in\{1,\ldots,k\}$, $s_i$ and
$s'_j$ are connected by the only flow $f_{ij}$ if $i<j$, by only
$f'_{ji}$ if $i>j$, and by the two flows $f_{ii}$ and $f'_{ii}$
if $i=j$. We show the following
  \begin{theorem}  \label{tm:sdecom}
(DS) can be solved in $O(E\log(2+V^2/E))$ time.
  \end{theorem}

This generalizes Theorem~\ref{tm:decom} for integer flows because a
digraph $D$ with an integer $S$--$T$ flow $g$ is turned into a
skew-symmetric graph with an integer symmetric
$(S\cup T')$--$(S'\cup T)$ flow by adding a disjoint copy of the
reverse to $D$ with the flow reverse to $g$ in it.

Our algorithm to solve (DS) relies on the following lemma (where, as
before, primes are used for the corresponding mate objects).
  \begin{lemma}  \label{lm:lem}
Let $g$ be a (not necessarily symmetric) half-integer flow from $S$
to $T$ in a skew-symmetric graph $G=(V,E)$ such that $\diver_g(v)$
is an integer for each $v\in V$. Let $g+g'$ be integer. Then there
exists, and can be found in $O(E)$ time, an integer flow $h$ in
$G$ such that $h+h'=g+g'$ and $\diver_h(v)=\diver_g(v)$ for all
$v\in V$.
  \end{lemma}
  \begin{proof}
Let $E_0$ be the set of arcs $e$ with $g(e)\not\in\Zset$. The
integrality of $g+g'$ implies $E'_0=E_0$, so the subgraph $\langle
E_0\rangle$ induced by $E_0$ is skew-symmetric. Also the
half-integrality of $g$ and the integrality of $\diver_g$ imply
that each node is incident with an even number of arcs in $E_0$. So
the underlying undirected graph $H$ of $\langle E_0\rangle$ is
Eulerian.

We grow a (simple) path $P$ in $\langle E_0\rangle$ such that
$P\cap P'=\emptyset$, starting with an arbitrary node $v_0$ and
allowing backward arcs in $P$.
Let $v$ be the last node of the current $P$,
and choose an arc $e\in E_0$ incident with $v$ and different from
the last arc of $P$ ($e$ exists as $H$ is Eulerian). Let $u$ be the
end of $e$ different from $v$. Three cases are possible. (i) If
both $u,u'$ are not in $P$, we increase $P$ by adding $e,u$, and
continue the process. (ii) If $u\in P$, we remove the part of $P$
from $u$ to $v$, obtaining the new current path from $v_0$ to $u$,
and add $e$ to the removed part, forming circuit $C$ (with possible
backward arcs). (iii) If $u'\in P$, we remove the part $Q$ of $P$
from $u'$ to $v$, obtaining the new current path, and add $e,Q'$
and $e'$ to $Q$, forming circuit $C$ (which is reverse to $C'$).

In case (ii), we update $g$ by pushing half-unit along $C$ (i.e.,
by setting $g(a):=g(a)+\frac12$ for the forward arcs $a$ in $C$,
and $g(a):=g(a)-\frac12$ for the backward arcs $a$) and by pushing
half-unit along the circuit reverse to $C'$. Accordingly, we update
$g'$ by pushing half-unit along $C'$ and along the circuit reverse
to $C$. And in case (iii), $g$ ($g'$) is updated by pushing
half-unit along $C$ (resp. $C'$). In both cases, the
new $g'$ is symmetric to the new $g$ and each of the functions
$g+g'$ and $\diver_g$ preserves. Also $E_0$ decreases by the set of
arcs occurring in $C\cup C'$, and the new $H$ is Eulerian. We
continue the process with the new $P$.

The final $g,g'$ give the desired $h,h'$. The bound $O(E)$ is
obvious.
  \end{proof}

Return to problem (DS). Add to $G$ new nodes $t,t'$ and arcs
$(t,s_i)$ and $(s'_i,t')$, forming skew-symmetric graph $G_1$, and
extend $f$ to an IS-flow from $t$ to $t'$ in $G_1$ in a natural way.
The fact that $f$ is integer and symmetric implies that $\diver_f(t)$
is even.

So we can apply Lemma~\ref{lm:lem} to the flows $g:=g':=\frac12 f$,
obtaining corresponding integer flows $h,h'$. The restriction
$\bar h$ of $h$ to $E$ is an integer flow from $S$ to $S'$, and we
apply the $O(E\log(2+V^2/E))$-algorithm from Section~\ref{sec:decom}
to decompose it as
  $$
  \bar h=\sum\nolimits_{1\le i,j\le k} h_{ij},
  $$
where $h_{ij}$ is an integer flow from $s_i$ to $s'_j$. Then the
flows $f_{ii}:=h_{ii}$ for $i=1,\ldots,k$, and $f_{ij}:=h_{ij}+
h'_{ji}$ for $1\le i<j\le k$ are as required, and
Theorem~\ref{tm:sdecom} follows.

\medskip
\noindent {\bf Remark.} The above proof involves the following
corollary from Lemma~\ref{lm:lem}.
  \begin{corollary} \label{cor:cor}
Let $f$ be an IS-flow from $S$ to $S'$ in a skew-symmetric graph $G
=(V,E)$ (where $|S|$ is not fixed), and let $\diver_f(v)$ be even
for all $v\in V$. Then there exists, and can be found in $O(E)$
time, an integer flow $g$ from $S$ to $S'$ such that $f=g+g'$ and
$\diver_f(v)=2\diver_g(v)$ for all $v\in V$.
  \end{corollary}


\end{document}